\documentclass{svproc}
\usepackage{url}

\usepackage{graphicx}

\newcount\eqnumber
\def\neweq{{\rm{(\the\eqnumber)}}\global\advance\eqnumber by 1}
\def\eqdef#1{\eqno\xdef#1{\the\eqnumber}\neweq}
\def\newaeq{{\rm{(\the\eqnumber a)}}\global\advance\eqnumber by 1}
\def\eqdaf#1{\eqno\xdef#1{\the\eqnumber}\newaeq}
\def\eqdisp#1{\xdef#1{\the\eqnumber}\neweq}
\def\eqdasp#1{\xdef#1{\the\eqnumber}\newaeq}

\newcount\refnumber
\def\newref{{\the\refnumber}\global\advance\refnumber by 1}
\def\refdef#1{{\xdef#1{\the\refnumber}}\newref}

\newcount\fignumber
\def\newfig{{\the\fignumber}\global\advance\fignumber by 1}
\def\figdef#1{{\xdef#1{\the\fignumber}}\newfig}

\def\smallskip{\vskip 3pt}
\def\medskip{\vskip 6pt}
\def\bigskip{\vskip 12pt}

\begin{document}
\mainmatter              
\title{On unconventional discretisations}
\titlerunning{Unconventional Discretisations}  
%
\author{Basil  Grammaticos\inst{1}, Thamizharasi Tamizhmani\inst{2} \and Ralph Willox\inst{3}}
\authorrunning{Grammaticos et al.} 
%
\tocauthor{Basil  Grammaticos, Thamizharasi Tamizhmani and Ralph Willox}
\institute{Universit\'e Paris-Saclay and Universit\'e Paris-Cit\'e, CNRS/IN2P3, IJCLab, 91405 Orsay, France\ \email{(bgrammat@ijclab.in2p3.fr)}
\and
SAS, Vellore Institute of Technology, Vellore - 632014, Tamil Nadu, India
\and 
Graduate School of Mathematical Sciences, the University of Tokyo, 3-8-1 Komaba, Meguro-ku, 153-8914 Tokyo, Japan }

\maketitle              

\begin{abstract}
We present various discretisation techniques that mark a departure from the classical algorithms of numerical analysis, in the aim of preserving the physical behaviour of the solutions to the differential systems we are discretising, even for discretisation steps that are not very small. These techniques are inspired by the principles laid down by Mickens and are close to those proposed by Hirota and Kahan. Detailed applications of these methods are presented for the logistic equation and for the Lotka-Volterra system. We discuss the importance of preserving positivity whenever the latter is an expected property of the solution and give a practical rule for achieving this. Possible problems arising from the discretisations we propose are also discussed. A section is devoted to discretisations of integrable systems aimed at preserving integrability. Finally we conclude this review with a discussion of the possible relevance of a discrete space-time.
\end{abstract}
\section{Introduction}
The mathematical description of physical laws is (almost) always based on differential equations. Understanding the principles underlying this formalism took, literally, centuries. In fact, five centuries before the current era, Zeno formulated the arrow (o\"{\i}stos) paradox which concludes that motion is impossible. Considering a flying arrow, Zeno states that, for the motion to occur, the arrow must change the position it occupies. But at every instant of time there is no motion occurring. So, if everything is motionless at every instant, and time is composed of instants, motion is indeed impossible. The fallacy in this reasoning lies in the assumption that time can be composed of instants of zero duration. Zeno tacitly assumes that the addition of infinitely many instants cannot be finite unless the instants have zero duration.
Our modern understanding of space-time is that a finite interval can be thought as the sum of intervals which can be as small as we like but not of strictly zero extent. The concept of limit, introduced in the infinitesimal calculus, then allows us to define instantaneous rates of variation.

Borrowing an example from Mechanics, if one wishes to study the motion of an object, one must solve a differential equation of the form
$${d^2{\vec r}\over dt^2}=F(\vec r, {d{\vec r}\over dt},t),\eqdef\zena$$
where $\vec r(t)$ is the position of the object at time $t$. In some, exceptional, cases the solution of (\zena) can be obtained analytically. However, this is not true for the majority of cases one encounters in the mathematical modelling of real-world situations. In that case one must resort to numerical simulations implemented on digital computers. The main problem with the latter is that while now-a-days they are lightning fast, grosso-modo, the only arithmetical operation they can carry out is addition. (The situation was different in the bygone era of analog computers which could solve equations involving derivatives, but had to be built specifically for a given problem, not at all like their multipurpose digital brethren).

Faced with the problem of simulating a dynamical equation on a digital machine, one must translate the differential system into a sequence of recursions involving simple arithmetical operations. This is where the notion of discretisation makes its appearance. The usual approach is to use ``black-box'' integrators, furnished by the standard libraries existing for the various computer languages that are in use today. The main drawback of this approach is that when one runs into problems, for example when the simulation results start differing wildly from the real physical behaviour, it is usually difficult to pinpoint their origin and implement adequate remedies. The aim of this paper is to show that, by following a few simple rules, one can develop quite efficient integrators, while keeping perfect control over the discretisation procedure. A further motivation of the present work can be attributed to the authors' interest in integrable systems. Given an integrable differential system, a question that is frequently asked is whether one can produce a discrete (difference) analogue that preserves the property of integrability. The answer to this question is often a resounding `yes', if one follows the guidelines we set out hereafter.

The paper opens with a presentation of a case study, that of the so-called `logistic' equation, for which we show that two classical discretisation schemes, Euler and Runge-Kutta, lead to unphysical results when the parameter of the equation exceeds some value, while a non-standard scheme proposed by Skellam and independently by Morishita, always performs in a satisfactory way. Before presenting any concrete discretisation techniques we first devote a section summarising the principles laid down by Mickens and which constitute a perfect guide for the discretisation of differential systems. Two other major approaches, those of Hirota and Kahan, are presented in the next section and are illustrated by several examples. The positivity condition is presented in the next section where we also discuss possible drawbacks of the methods introduced in this paper. 
We claim that the unconventional discretisation methods, introduced here, lead to physical results even when the discretisation step is large. In fact, the limit where the discretisation step goes to infinity has meaning. It is linked to the so-called `ultradiscretisation' procedure, which allows one to construct a celllular-automaton analogue of a given discrete system. 
A final section is devoted to the question of integrable discretisations and we conclude with a discussion on the purported continuous or discrete character of space-time.
\section{A case study}
The first case we are going to study is the discretisation of the logistic equation, an equation first proposed by Verhulst [\refdef\verhu] in the middle of the 19th century as a model for population growth under the constraint of scarcity of natural resources. It has the form of a Riccati equation
$$x'=ax-bx^2,\eqdef\zdyo$$
where $a$ and $b$ are positive constants, and its analytical solution can be easily obtained (introducing the variable $1/x$ transforms the equation into a linear one). We find
$$x={ax_0\over bx_0+(a-bx_0)e^{-a(t-t_0)}}.\eqdef\ztri$$
Starting form a value $x_0<a/b$ at $t=t_0$, $x$ increases monotonically towards the value $a/b$. 

A first discretisation we are going to consider is a simple Euler scheme using a forward difference for the first derivative. We have
$${x_{n+1}-x_n\over \delta}=ax_n-bx_n^2,\eqdef\ztes$$
where $\delta>0$ is the time-discretisation step. We introduce a new parameter $r=1+a\delta>1$ and scale $x$ as $x\mapsto r x / (b\delta)$ to bring the equation to the form 
$$x_{n+1}=rx_n(1-x_n).\eqdef\tdyo$$
This is the well-known logistic map which has been shown to behave chaotically for values of the parameter $r$ sufficiently close to (but not greater than) $r=4$ [\refdef\logis, \refdef\may]. When $r>4$ the map  (\tdyo) can generate values greater than 1, even for initial values $0<x_0<1$, and the positivity of successive iterates breaks down.

A second discretisation will be based on a second-order Runge-Kutta scheme [\refdef\abram], often used in numerical simulations. Given a first-order differential equation of the form $x'=f(x)$ this scheme proposes the discretisation $x_{n+1}=x_n+\delta(f(x_n)+f(x_n+\delta f(x_n))/2$. Applying the scheme to the logistic equation with the choices $r=1+a\delta$ and the above scaling (allowing us to take $b=r/\delta$), we find for the Runge-Kutta scheme the recursion
$$x_{n+1}={1\over2}\big(x_n+r^2x_n(1-x_n)-r^3x_n^2(1-x_n)^2\big).\eqdef\zhex$$
However a third, unconventional, discretisation due to Skellam [\refdef\skel] and, independently, Morishita [\refdef\moris] does also exist. It consists in replacing the $x_n^2$ term by $x_nx_{n+1}$, which at the continuum limit goes to $x^2$ just as $x_n^2$ does. (Later in the article we will discuss the possible consequences of such a choice). Using again a forward difference for the derivative and the above scaling of $x$, we find the recursion
$$x_{n+1}={rx_n\over 1+rx_n}.\eqdef\zhep$$
Before proceeding to numerical simulations of the above schemes it is interesting to study their fixed points. Besides a trivial unstable fixed point at $x=0$, both the logistic mapping and the one due to Morishita and Skellam have a unique non-zero fixed point at $x_*=1-1/r$, for any value of $r>1$. For the logistic mapping this fixed point is attractive as long as $r<3$, but it becomes unstable for $r>3$. For the mapping (\zhep), on the other hand, this fixed point is attractive for any value $r>1$ exactly as the corresponding fixed point for the logistic equation (\zdyo). The Runge-Kutta scheme has, besides a trivial (repulsive) fixed point at $x=0$, a unique attractive fixed point at $1-1/r$ when $r<3$. This fixed point, however, becomes repulsive when $r>3$ at which value a typical pitchfork-style bifurcation occurs and two extra fixed points appear. These are given by the roots of the equation $r^2x_*^2-r(r+1)x_*+r+1=0$ which are both positive when $r>3$ (and complex when $r<3$). Both these new fixed points remain attractive as long as $3<r<1+2\sqrt{2}$ but become unstable when $r>1+2\sqrt{2}$. 

We proceed now to simulations using the schemes introduced above. In Figure \figdef\one\ we show the results for $r=2$.
\medskip
\centerline{\includegraphics[width=10 cm,keepaspectratio]{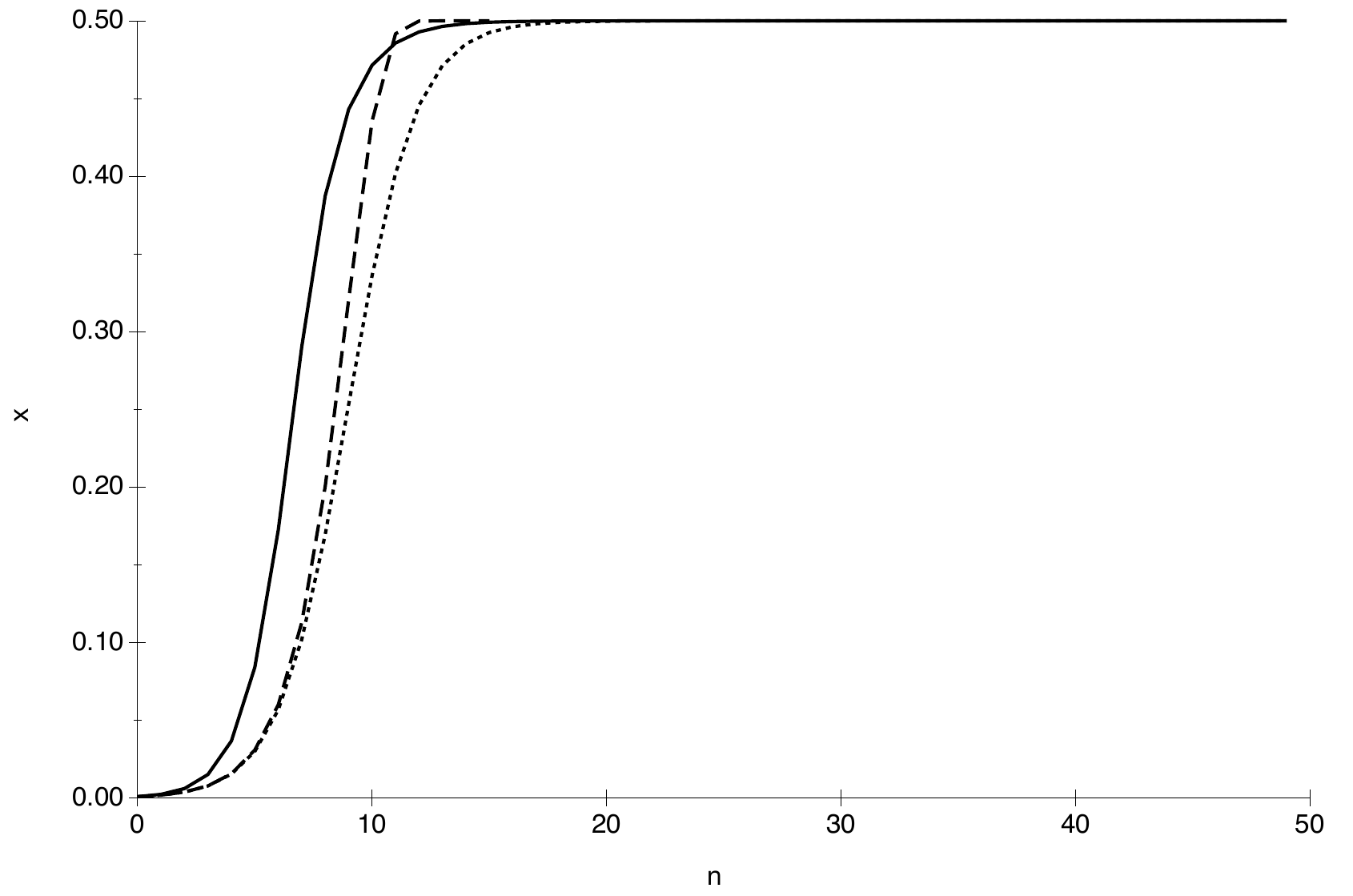}}
\smallskip{\bf Figure \one}. {\sl Numerical solution of the logistic equation, for $r=2$, using the Euler scheme (\tdyo) [dashed line], the Runge-Kutta scheme (\zhex) [continuous line] and the Morishita scheme (\zhep) [dotted line]}.
\smallskip
We remark that all three schemes represent the dynamics of the logistic equation adequately. Next we turn to the results for $r=3$ and represent them graphically in Figure \figdef\two.
\medskip
\centerline{\includegraphics[width=10 cm,keepaspectratio]{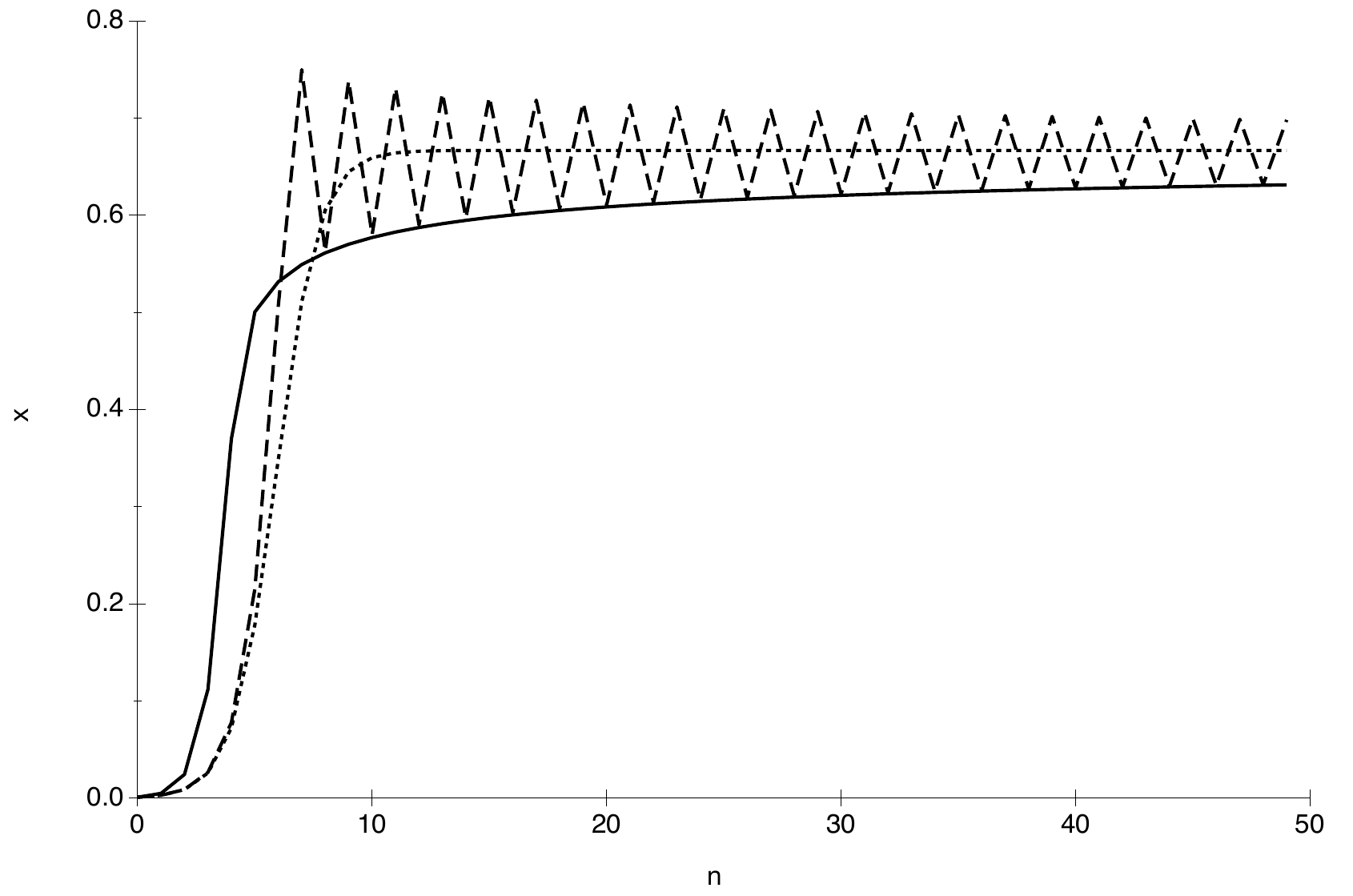}}
\smallskip{\bf Figure \two}. {\sl Numerical solution of the logistic equation, for $r=3$, using the Euler scheme (\tdyo) [dashed line], the Runge-Kutta scheme (\zhex) [continuous line] and the Morishita scheme (\zhep) [dotted line]}.
\smallskip
The solution of the Euler scheme has now entered an oscillation of period 2, while for the Runge-Kutta scheme the fixed point $1-1/r$ is about to lose its stability. The Morishita scheme however still provides a satisfactory description of the dynamics of the logistic equation. As is well-known, for $r>3$ the logistic mapping (Euler scheme) enters in a period doubling cascade which leads to a chaotic behaviour and, as mentioned above, for $r>4$ its solution ceases to be positive. In Figure \figdef\three\ below we show the results of the numerical simulation for $r=4$.
\medskip
\centerline{\includegraphics[width=9 cm,keepaspectratio]{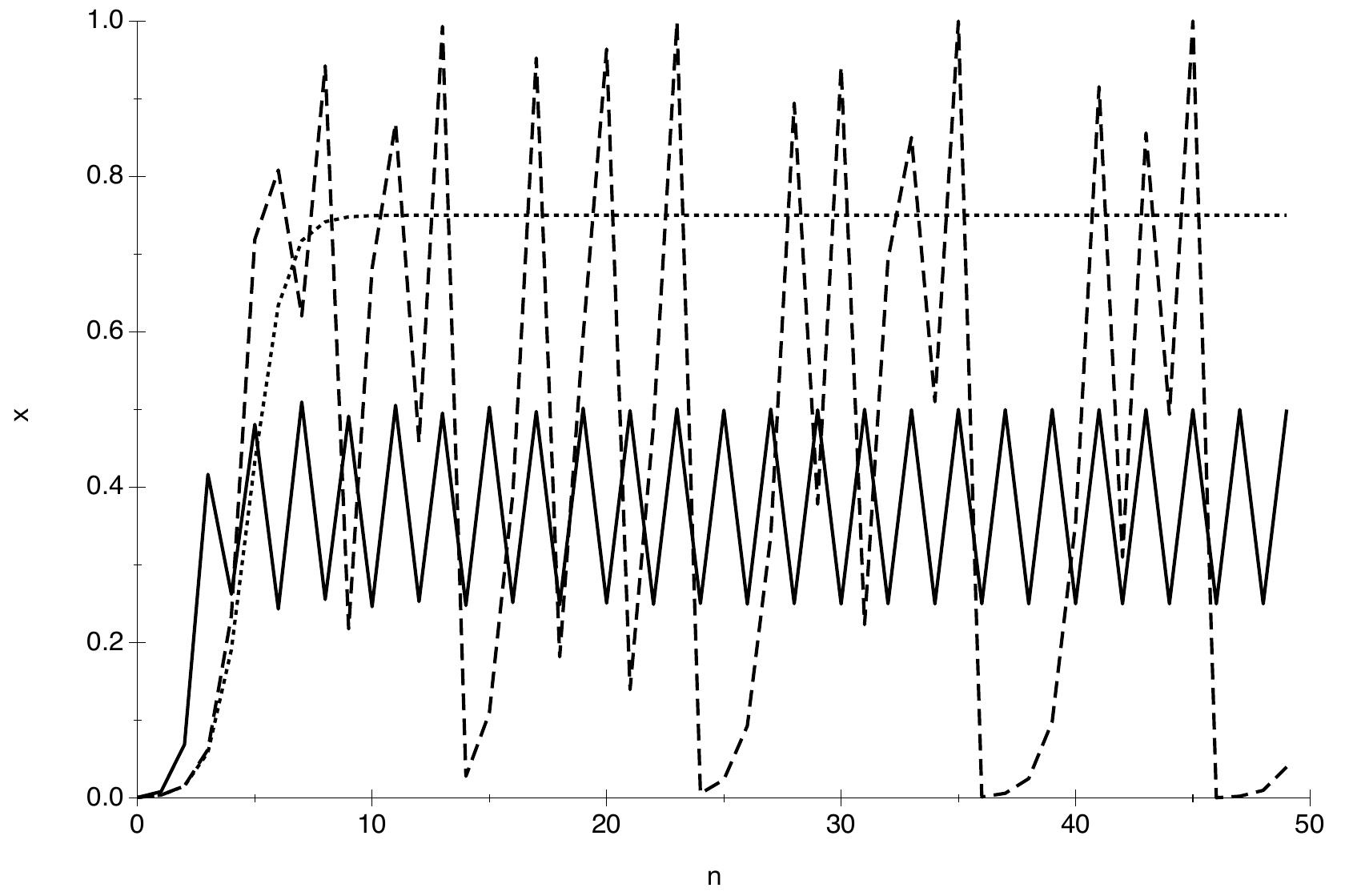}}
\smallskip{\bf Figure \three}. {\sl Numerical solution of the logistic equation, for $r=4$, using the Euler scheme (\tdyo) [dashed line], the Runge-Kutta scheme (\zhex) [continuous line] and the Morishita scheme (\zhep) [dotted line]}.
\smallskip
We remark that the solution for the Runge-Kutta scheme is no more attracted to the fixed point  $1-1/r$  but is now oscillating around the smaller of the two  new fixed points. In fact, the fixed point $1-1/r$ acts as a (0-dimensional) separatrix: when $r>3$, iterates of initial conditions that are greater than $1-1/r$ are attracted to the larger stable fixed point and those for initial conditions that are smaller than  $1-1/r$ to the smaller of the two new fixed points.When $r<1+2\sqrt{2}$ the iterates converge to said fixed points but when $r>1+2\sqrt{2}$ they exhibit a period-2 oscillation around them, signalling a possible route to chaos. In fact, for $r$ greater than a value around 4.4037, the solution provided by the Runge-Kutta scheme is no more positive but more interestingly, for values of $r$ slightly smaller than this,  the behaviour of the mapping appears to be chaotic with an interesting interplay between the two fixed points. In Figure \figdef\four\ we show how the solution, starting from a point close to the larger of the two fixed points is finally attracted to the lower of the two. 

\medskip
\centerline{\includegraphics[width=10 cm,keepaspectratio]{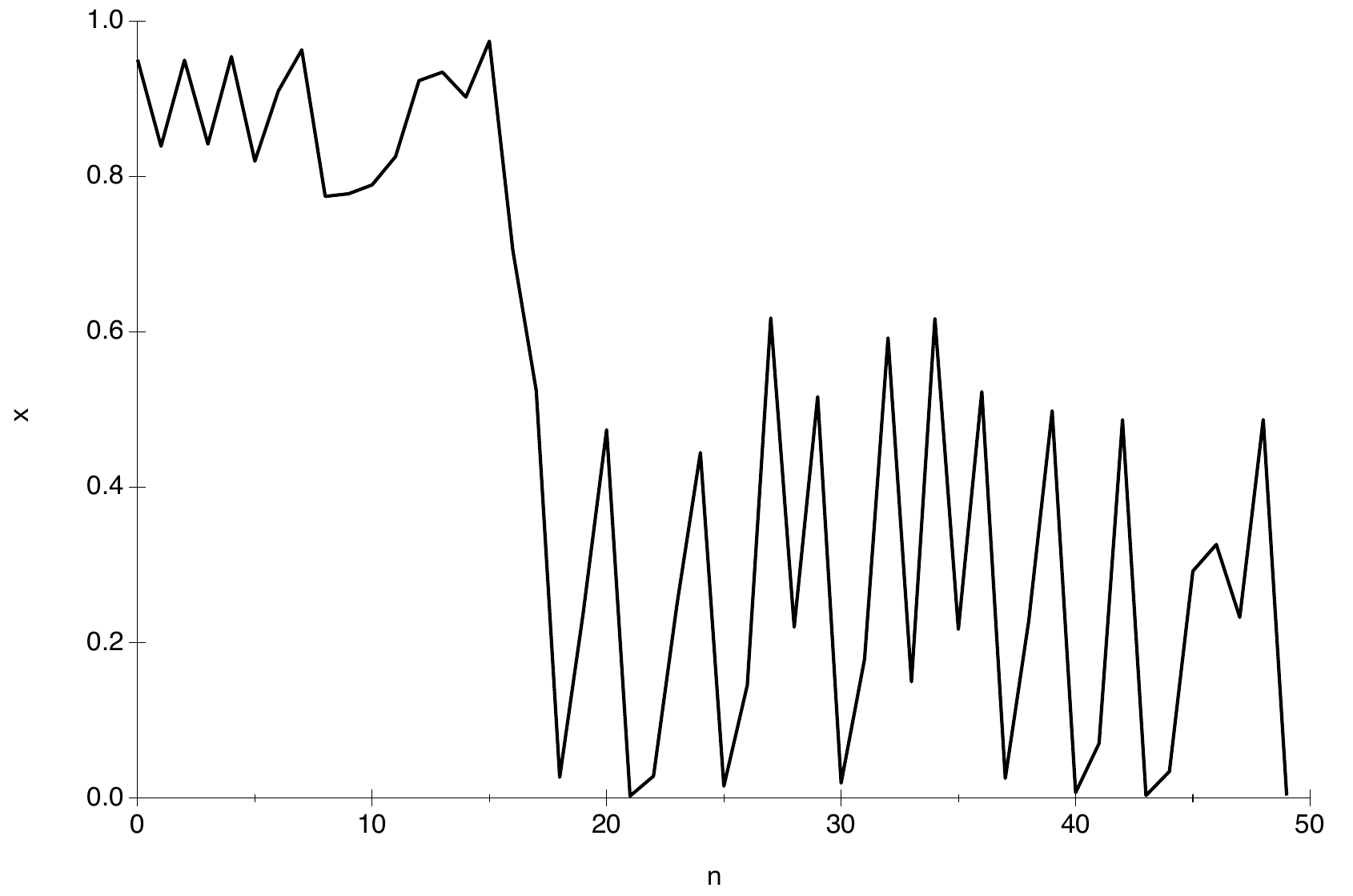}}
\smallskip{\bf Figure \four}. {\sl Numerical solution of the logistic equation, for $r=4.4$, using  the Runge-Kutta scheme}.
\smallskip
It goes without saying that for values of $r$ larger than the one where the Runge-Kutta scheme ceases to furnish positive solutions, the Skellam-Morishita scheme (\zhep) still, unfailingly, provides a monotonous solution that converges asymptotically towards the fixed point $1-1/r$.

We are not going to delve further into the behaviour of the ``classical'' integration schemes, Euler and Runge-Kutta. (The behaviour of the latter would justify further studies but this is, alas, outside of the scope of this paper). They both work fine for small values of the parameter $r$ (which we remind is directly related to the discretisation step $\delta$), close to 1, but lead to an ``unphysical'' behaviour for larger values of said parameter. On the other hand the unconventional discretisation of Skellam-Morishita always works, guaranteeing the positivity of  the solution ( but more on this later).
\section{The Mickens discretisation principles}
The examples of the previous section show that trusting standard integrators is not necessarily the optimal choice. On the other hand devising discretisation methods ``out of the blue'' without at least some general principles is not something that is particularly recommendable either. Fortunately when it comes to discretisations, R. Mickens [\refdef\ronal]  has not only pioneered the approach, introducing methods that have definite implications in the domain of numerical analysis, but he has summarised his experience in a series of simple principles (``{\sl distilling the essence of the wisdom discovered thus far}', in his own words).
The principles of Mickens should serve as a guide in any attempt at discretising differential systems, and all the more so if the approach is along the lines we dubbed ``unconventional''. We present them below, where in italics we quote Mickens {\sl verbatim} [\refdef\micke] and follow each principle with a few comments of our own. 

{\sl i) The orders of the discrete derivatives should be equal to the orders of the corresponding derivatives appearing in the differential equations.}

In general, when the order of the difference equation is larger than the order of the differential equation, numerical instabilities that may oscillate with increasing amplitude will appear. Among the standard methods proposed for the integration of differential equations one finds, for example in the Abramowitz-Stegun monograph [\abram], a discretisation of the first order equation $x'=f(x)$ through a supposedly improved Euler method where the recursion $x_{n+1}=x_n+\delta f(x_n)$, accurate at order $\delta^2$, is replaced by $x_{n+1}=x_{n-1}+2\delta f(x_n)$, in principle accurate at order $\delta^3$. This is the kind of ``improvement'' that Mickens is warning against, since the order of the difference equation is not the same as that of the differential one.

{\sl ii) Discrete representations for derivatives must, in general, have nontrivial denominator functions.}

Usually when it comes to discretising the first derivative $x'$ one uses a simple forward difference $(x_{n+1}-x_n)/\delta$. What Mickens recommends here is that, in some cases one should think about a denominator that is a function of $\delta$, i.e. introduce a discretisation $(x_{n+1}-x_n)/d(\delta)$ where for $\delta$ sufficiently small $d(\delta)$ is of order $\delta$. In the next section  we shall encounter such an example.
However, as Mickens himself points out, the particular form of the denominator function is not a critical feature.

{\sl iii) Nonlinear terms should, in general, be replaced by nonlocal discrete representations.}

For instance, terms of the form $x^2$ and $x^3$ which may occur in a first-order, nonlinear ordinary differential equation could be modelled as follows: $x^2\to x_nx_{n+1}$ and $x^3\to x_n^2(x_n+x_{n+1})/2$. As we have seen in the previous section, the choice of Morishita using $x_nx_{n+1}$ instead of $x_n^2$ ensures the stability to his integrator. Moreover the application of this principle of Mickens has the advantage, in the present case, to guarantee the positivity of the solution. As we shall see later in this paper positivity is something that can, in the proper context, be elevated to the status of a principle.

{\sl iv) All properties of the continuous equation must be present in the discrete one.}

To be fair we are slightly getting ahead of ourselves on this point. Mickens presents his fourth principle as {\sl Special conditions that hold for either the differential equation and/or its solutions should also hold for the difference equation model and/or its solutions} and offers the example of time-reversal invariance which, if present in the differential system, should also be present in the discrete one. However he immediately forges ahead to introduce the principle of ``dynamic consistency'' which adopts the broad formulation of iv). Mickens gives a (non-exhaustive) list of  properties that should be preserved by the discretisation like positivity,  boundedness,  monotonicity, fixed-points and their stability properties, integer valued dependent variables, existence of special solutions, limit cycles etc. Given that the principle of dynamic consistency can, when adequately interpreted, encompass the more technical rule concerning special conditions, we deem it justified to just focus on the former. Mickens accompanies his principles with a prudent caveat. In his own words: ``{\sl In general, one should not expect the above stated rules to give a best difference equation model for any particular differential equation of interest. However, we do expect that their application will allow the elimination of many instabilities that can occur if the usual rules are followed for constructing finite-difference schemes for differential equations}''. 

When one is given a differential equation there exist a profusion of discrete schemes that obey the first three principles. However when one asks for the properties of the differential system to carry over to its discrete analogue the available choice is pared down considerably (and sometimes none survives). One domain of particular interest to the authors is that of integrability. Applied to this particular case the principle dictates that the discrete analogue of a continuous integrable system be integrable as well. While this is a tall order, pursuing this line of research leads to very interesting results.
\section{The Hirota and the Kahan discretisations}
Hirota is famous for his construction of $N$-soliton solutions [\refdef\hirot] for a collection of integrable two-dimensional evolution equations. To this end he introduced the bilinear formalism that constitutes a precious tool for the study of integrable systems. What is less known is that Hirota, from the very outset of his research in integrable systems, also derived discrete versions of said evolution equations based on the requirement that the $N$-soliton solution be preserved by the discretisation. (One sees immediately that this is an application of the 4th principle of Mickens, and in fact a rather exacting version thereof).

The implications of requiring that the discretisation preserves the solution of the differential system can be nicely illustrated through the example of the logistic equation. The construction that follows is not due to Hirota. It is the one proposed by Morishita [\moris], leading to the mapping (\zhep). From equation (\ztri) we have the exact solution of the logistic equation and proceed to discretise time by putting $t=t_0+n\delta$. We find
$$x_n={cx_0\over x_0+(c-x_0)e^{-na\delta}},\eqdaf\zoct$$
with $c=a/b$, and write the solution at the next time step
$$x_{n+1}={cx_0\over x_0+(c-x_0)e^{-(n+1)a\delta}}.\eqno(\zoct {\rm b})$$
Next we eliminate the explicit dependence on $n$, i.e. the term $e^{-na\delta}$, between the two equations and obtain the discrete equation 
$${x_{n+1}-x_n\over1-e^{-a\delta}}=(1-{x_n\over c})x_{n+1},\eqdef\zenn$$
which simplifies to
$$x_{n+1}= {e^{a\delta}x_n\over1+ {(e^{a\delta}-1)\over c} x_n}.\eqdef\zenm$$
A remark is in order at this point. The denominator of the discrete derivative appearing in the left hand side of (\zenn) is a perfect example of the third Mickens principle. If one wishes to preserve the solution for any choice of the values of the parameters and the discretisation step, one must keep the full denominator as appearing in (\zenn). Making the assumption that the product $a\delta$ is small enough, one can replace $1-e^{-a\delta}$ by simply $a\delta$ and (after scaling $x$) use the fact that $c=(r-1)/r$ to obtain mapping (\zhep).

In order to illustrate Hirota's bilinear approach we choose a very simple example, that of a general Riccati equation.  We start from 
$$x'=ax^2+2dx+f,\eqdef\zdek$$
and introduce the ansatz $x=P/Q$.
We remark that a gauge transformation $P\to \phi P$, $Q\to \phi Q$ leaves $x$ invariant.
Substituting the ansatz into (\zdek) we find
$$PQ'-QP'=aP^2+2dPQ+fQ^2,\eqdef\dena$$
which is the so-called `bilinear form' of (\zdek), and is indeed invariant under the above gauge transformation.
In order to discretise (\dena) we introduce a discrete time step $\delta$ and propose
$${Q_{n+1}P_n-P_{n+1}Q_n\over\delta}=aP_n^2+2dP_nQ_n+fQ_n^2.\eqdef\ddyo$$
However (\ddyo) is not gauge-invariant. Hirota restores this invariance by introducing a special discretisation of the quadratic terms leading to:
$${Q_{n+1}P_n-P_{n+1}Q_n\over\delta}=aP_nP_{n+1}+d(\alpha Q_{n+1}P_n+\beta P_{n+1}Q_n)+fQ_nQ_{n+1}\eqdef\dtri$$
where $\alpha+\beta=2$. Dividing by $Q_nQ_{n+1}$ and introducing the variable $x_n=P_n/Q_n$ we finally obtain
$${x_{n+1}-x_n\over\delta}=ax_nx_{n+1}+d(\alpha x_n+\beta x_{n+1})+f.\eqdef\dtes$$
Hirota's gauge-invariant discretisation can be applied to any bilinear system. However one should keep in mind that, given a nonlinear system, it is not always easy to obtain the proper bilinearisation.

Kahan's method was first presented in unpublished notes [\refdef\kahan] under the moniker of ``unconventional'' and went mostly unnoticed outside the domain of numerical analysis, to be re-discovered rather recently by the integrability community. It concerns the discretisation of quadratic nonlinearities. We will illustrate it using again the example of Riccati equation. Starting from (\zdek) Kahan proposes the discrete form
$${x_{n+1}-x_n\over\delta}=aq(x_n,x_{n+1})+d(x_n+x_{n+1})+f,\eqdef\dpen$$
where the quadratic form $q(x_n,x_{n+1})$, in this one-dimensional case, is simply
$$q(x_n,x_{n+1})={1\over2}\big( (x_n+x_{n+1})^2-x_n^2-x_{n+1}^2\big)=x_nx_{n+1}.\eqdef\dhex$$
We remark that the resulting discrete equation is quite similar to the one obtained by Hirota. 

It is interesting to confront the two methods in the case of the Lotka-Volterra predator-prey system
$$x'=x(\lambda-y)\eqdaf\dhep$$
$$y'=y(x-\mu),\eqno(\dhep{\rm b})$$ 
where $\lambda, \mu>0$.
The discretisation following Kahan's prescription leads to
$${x_{n+1}-x_n\over\delta}=-(x_ny_{n+1}+y_nx_{n+1})/2+\lambda(x_n+x_{n+1})/2\eqdaf\doct$$
$${y_{n+1}-y_n\over\delta}=(x_ny_{n+1}+y_nx_{n+1})/2-\mu(y_n+y_{n+1})/2.\eqno(\doct{\rm b})$$ 
Solving for $x_{n+1}$ and $y_{n+1}$ and introducing $\epsilon=\delta/2$ we find
$${x_{n+1}\over x_n}={x_n\epsilon(1+\epsilon\lambda)+y_n\epsilon(1-\epsilon\mu)-(1+\epsilon\lambda)(1+\epsilon\mu)\over
x_n\epsilon(1-\epsilon\lambda)-y_n\epsilon(1+\epsilon\mu)-(1-\epsilon\lambda)(1+\epsilon\mu)}\eqdaf\denn$$
$${y_{n+1}\over y_n}={-x_n\epsilon(1+\epsilon\lambda)-y_n\epsilon(1-\epsilon\mu)-(1-\epsilon\lambda)(1-\epsilon\mu)\over
x_n\epsilon(1-\epsilon\lambda)-y_n\epsilon(1+\epsilon\mu)-(1-\epsilon\lambda)(1+\epsilon\mu)}.\eqno(\denn{\rm b})$$
A relative drawback of equation (\denn) is that when the discretisation step is large one may obtain negative values for $x$ and $y$. 

The Hirota discretisation based on the bilinear formalism, together with the gauge-invariance requirement gives the discrete system
$${x_{n+1}\over x_n}={x_n\delta\beta(1+\delta\lambda_2)+y_n\delta\beta(1-\delta\mu_2)-(1+\delta\lambda_2)(1+\delta\mu_1)\over x_n\delta\beta(1-\delta\lambda_1)-y_n\delta\alpha(1+\delta\mu_1)-(1-\delta\lambda_1)(1+\delta\mu_1)}\eqdaf\vdek$$
$${y_{n+1}\over y_n}={x_n\delta\alpha(1+\delta\lambda_2)+y_n\delta\alpha(1-\delta\mu_2)+(1-\delta\lambda_1)(1-\delta\mu_2)\over -x_n\delta\beta(1-\delta\lambda_1)+y_n\delta\alpha(1+\delta\mu_1)+(1-\delta\lambda_1)(1+\delta\mu_1)}.\eqno(\vdek{\rm b})$$
where $\alpha+\beta=1$, $\lambda_1+\lambda_2=\lambda$ and $\mu_1+\mu_2=\mu$. Hirota's  discretisation has more freedom  compared to Kahan's and can, with the appropriate choice of the parameters, lead to positive definite $x_n$ and $y_n$. But to be fair, Hirota's scheme is somewhat overburdened by parameters. It is thus interesting to introduce another ``unconventional''  method of discretisation that can lead to a much simpler and still perfectly satisfactory discrete scheme [\refdef\murat].
We start from (\dhep a) and in order to guarantee the positivity we discretise $xy$ as to $x_{n+1}y_n$ while the $\lambda x$ term becomes $\lambda x_n$. We find thus
$${x_{n+1}\over x_n}={1+\delta\lambda\over1+\delta y_n}.\eqdaf\vena$$
Similarly for (\dhep b) the positivity requirement suggests that the term $-\mu y$ become $-\mu y_{n+1}$. However two choices appear possible for the $xy$ term, either $x_ny_n$ or $x_{n+1}y_n$. It turns out that only the second guarantees the reversibility of the system, i.e. a backwards evolution expressed as a rational mapping. We have finally
$${y_{n+1}\over y_n}={1+\delta x_{n+1}\over1+\delta\mu}.\eqno(\vena{\rm b})$$
With these choices, the evolution of (\vena) always leads to positive values for $x$ and $y$ (provided we start from positive initial conditions and parameters, of course). Our scheme  is a special case of Hirota's, obtained from (\vdek) by taking $\alpha=1$, $\beta=0$, $\lambda_1=0$, $\lambda_2=\lambda$, $\mu_1=\mu$ and $\mu_2=0$. 

Another model that we will examine with the help of the Hirota and the Kahan approaches is that of the quartic oscillator for which the aim is to obtain a discretisation that has a conserved quantity, just as the original continuous system. 
Starting from the Hamiltonian 
$$H={1\over2}x'^2+{c\over4}x^4,\eqdef\vdyo$$ 
one obtains the equation of motion
$$x''+cx^3=0.\eqdef\vtri$$
Hirota presents two different approaches for the discretisation of (\vtri). The first is based on an a discretisation that ensures the conservation of the Hamiltonian. The following form is proposed
$$H={1\over2\delta^2}(x_n-x_{n-1})^2+{c\over4}x_n^2x_{n-1}^2,\eqdef\vtes$$
leading to the equations of motion
$${1\over\delta^2}(x_{n+1}-2x_n+x_{n-1})+{c\over2}x_n^2(x_{n+1}+x_{n-1})=0.\eqdef\vpen$$
This discretisation is in fact one proposed initially by Potts [\refdef\potts], although in their monograph [\refdef\takah], Hirota and Takahashi attribute it to Suris [\refdef\suris]. The second discretisation is based on the bilinear approach of Hirota. We shall not present all the calculational details but give directly the result. The discrete equation of motion is now
$${1\over\delta^2}(x_{n+1}-2x_n+x_{n-1})+cx_{n+1}x_nx_{n-1}=0.\eqdef\vhex$$
We introduce a scaling for $x$ allowing to put $c\delta^2=-1$, and find that (\vhex) has the conserved quantity
$$H={x_n^2x_{n-1}^2-2(x_n^2+x_{n-1}^2)+4\over x_nx_{n-1}-1}.\eqdef\vhep$$
This suggests a canonical, albeit less intuitive, form for the equation of motion
$$(x_nx_{n+1}-1)(x_nx_{n-1}-1)=1-2x_n^2,\eqdef\voct$$
which belongs to Class IV of Quispel-Roberts-Thompson (QRT) mappings [\refdef\qrt] in the classification given in [\refdef\class]. 

The application of the Kahan discretisation to the quartic oscillator runs into the problem that, originally, Kahan proposed his method only for the case of quadratic nonlinearities. Fortunately the question of higher nonlinearities was addressed by Celledoni and collaborators [\refdef\polar] and they proposed a solution through the introduction of what they call ``polarisation''. The third order polarisation for a given cubic form $C(x)$ is through the symmetric trilinear form 
$$\tilde C(x,y,z)={1\over 6}\big(C(x+y+z)-C(x+y)-C(y+z)-C(z+x)+C(x)+C(y)+C(z)\big),\eqdef\venn$$
where we used the tilde in order to denote the ``polarisation'' of the cubic form. 

In [\refdef\yuris] Suris gives a direct application of the Kahan cum polarisation method to the case of the quartic oscillator. The discrete equation of motion is identical to the one obtained by the Hirota method, equation (\vhex). And, using the latter, Suris finds that the conserved quantity can be written in a very simple form
$$H={1\over\delta^2}\big(x_n(x_{n+1}+x_{n-1})-2x_{n+1}x_{n-1}\big)\eqdef\vdek$$
involving the values of $x$ at three adjacent points.

We do not claim that  there are no other discretisation methods that will produce discrete systems that have conserved quantities similar to the original continuous model they are supposed to describe, but we do wish to emphasize that the methods we presented here are specifically designed so as to produce discrete models that do have such properties.
\section{Tricks and caveats}
In the section devoted to the logistic equation we saw that the scheme due to Skellam and Morishita not only guarantees a positive solution but also that it displays the (hoped for) monotonous behaviour towards an asymptotically constant value even for large values of the parameter $r$.  Similarly, in the previous section we presented our unconventional discretisation of the Lotka-Volterra equation, taking particular care to preserve the positivity of the solution. Thus we guarantee that the solutions of both discrete systems (logistic and Lotka-Volterra) will remain `physical' since they have to do with population densities. 

We feel that the positivity property is a very strong one which must be enforced whenever this is justified. In [\refdef\handy] we have formulated this in the form of an aphorism: ``if the dependent variables of an equation correspond to positive quantities, the discretisation should be such as to avoid any minus signs''. However, in that paper, which deals with the dynamics of the interaction between a homogeneous population and the natural resources and reserves needed for its survival, we ran into a difficulty when trying to apply this principle.  One, very, simple model we proposed in [\handy] was given by the differential system
$${dx\over dt}=\alpha z-\beta x\eqdaf\tena$$
$${dy\over dt}=y(1-y)-xy\eqno(\tena{\rm b})$$
$${dz\over dt}=xy-\gamma x.\eqno(\tena{\rm c})$$
which, in view of numerical simulations, we proceeded to discretise. Note that $\alpha, \beta$ and $\gamma$ are positive constants. In order to avoid any minus signs we staggered the first two equations as:
$${x_{n+1}-x_n\over\delta}=\alpha z_n-\beta x_{n+1}\eqdaf\tdyo$$
$${y_{n+1}-y_n\over\delta}=y_n-y_ny_{n+1}-x_ny_{n+1}.\eqno(\tdyo b)$$
However, no staggering can bring equation (\tena c) to a discrete form without minus signs. The trick to deal with this difficulty is to multiply the $\gamma x$ term by a fictitious factor $z/z$ and then apply the adequate staggering, as in:
$${z_{n+1}-z_n\over\delta}=x_ny_{n+1}-\gamma x_n {z_{n+1}\over z_n}.\eqno(\tdyo c)$$
The $xy$ term was discretised to $x_ny_{n+1}$ so as to coincide with the term present in (\tdyo b) but we could have equally taken $x_ny_n$, without this affecting the results. With the above staggering no minus signs appear in the ensuing discrete evolution equations. Thus if we start from positive initial conditions and positive parameters, it is guaranteed that the iteration of (\tdyo) will yield positive values at all times. 
Notice that, in order to ensure positivity through the introduction of the $z_{n+1}/z_n$ term, the discrete equations cease to be time-reversible, while the initial differential system did possess that property.  We believe that is a small price to pay, since our main interest will lie in the future evolution of the system.

At this point a caveat is in order. In all cases presented up to now we did not hesitate to consider the discrete variable corresponding to $x$ at either point $n$ or $n+1$. But one should not forget that, in the discrete time evolution, a contribution $x_n$ is effectively `delayed' with respect to $x_{n+1}$ by one time step. Hence, a discretisation that allows, as per Mickens' principles, for a non-local representation of the nonlinear terms, actually introduces delays which may affect the dynamical behaviour of the system. In [\refdef\delay] two of the present authors examined systems where these ``hidden'' delays do indeed affect the dynamics. And what makes the situation even worse is that when one introduces discretisations so as to do away with the hidden delays, the resulting system may have instabilities that upset the dynamical behaviour. These situations are admittedly special, but the caveat is quite justified and particular care is always necessary when one tackles the problem of discretisation. 

Coming back to the question of positivity, we have seen in the case of the logistic equation how to go from the logistic mapping to that of Skellam and Morishita: it suffices to discretise the quadratic term $x^2$ as $x_nx_{n+1}$ instead of $x_n^2$. So the natural question is whether something analogous can be done in the case of the Runge-Kutta scheme. The answer is obviously ``yes''. A look at the right hand side of (\zhex) reveals the presence of negative terms $x_n^2$ and $x_n^4$. It suffices then to introduce the discretisations $x_nx_{n+1}$ and $x_n^3x_{n+1}$ in order to transform (\zhex) into 
$$x_{n+1}=x_n{1+r^2+2r^3x_n^2\over2+r^2(r+1)x_n+r^3x_n^3}.\eqdef\ttri$$
Given the form of (\ttri) it is guaranteed that the solution will be positive at all time steps for any value of the parameter $r$. In figure \figdef\five\  we present the results of the iteration of (\ttri) starting from a small value of $x$ just as in the figures \one\ to \three\ for $r=4.4$.

We remark that, while the solution obtained by the ``standard'' Runge-Kutta is on the brink of becoming negative, the solution obtained by the scheme given by equation (\ttri) is positive and monotonous. However all the drawbacks of the Runge-Kutta scheme cannot be addressed by just enforcing positivity: the fixed point $x_*=1-1/r$ is still repulsive and the solution of (\ttri) is attracted to $x_*=0.3012$, i.e. the smaller root of the quadratic equation $r^2x_*^2-r(r+1)x_*+r+1=0$ for $r=4.4$. 
\medskip
\centerline{\includegraphics[width=8 cm,keepaspectratio]{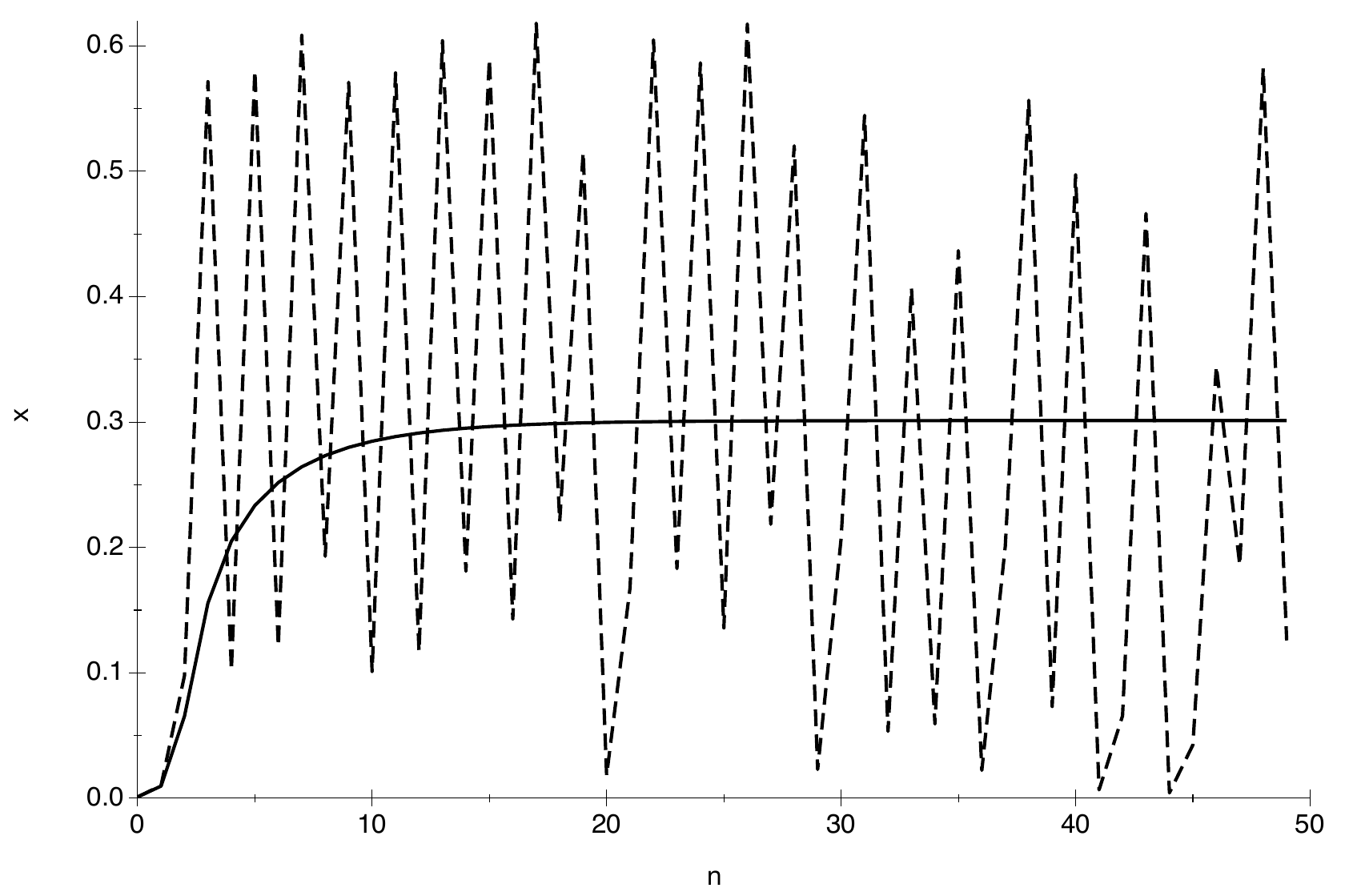}}
\smallskip{\bf Figure \five}. {\sl Numerical solution of the logistic equation, for $r=4.4$, using the Runge-Kutta scheme (\zhex) [dash-dot line] and the positive Runge-Kutta (\ttri) [continuous line]}.
\smallskip

As pointed out in section 2, the parameter $r$ is related to the discretisation step $\delta$ and increasing $r$ can be interpreted as an increase of the step $\delta$ (and a change of the scaling of the dependent variable). It is thus interesting to investigate the effect of the choice of time step on the dynamics. To this end we consider the discrete Lotka-Volterra scheme (\vena) introduced in the previous section. We iterate the discrete equations starting from the same initial conditions and with the same parameters $\lambda=1$ and $\mu=2$ for two choices of the time step: $\delta=0.001$ and $\delta=1$. Figure \figdef\six\ shows the results of this simulation. 
\medskip
\centerline{\includegraphics[width=10 cm,keepaspectratio]{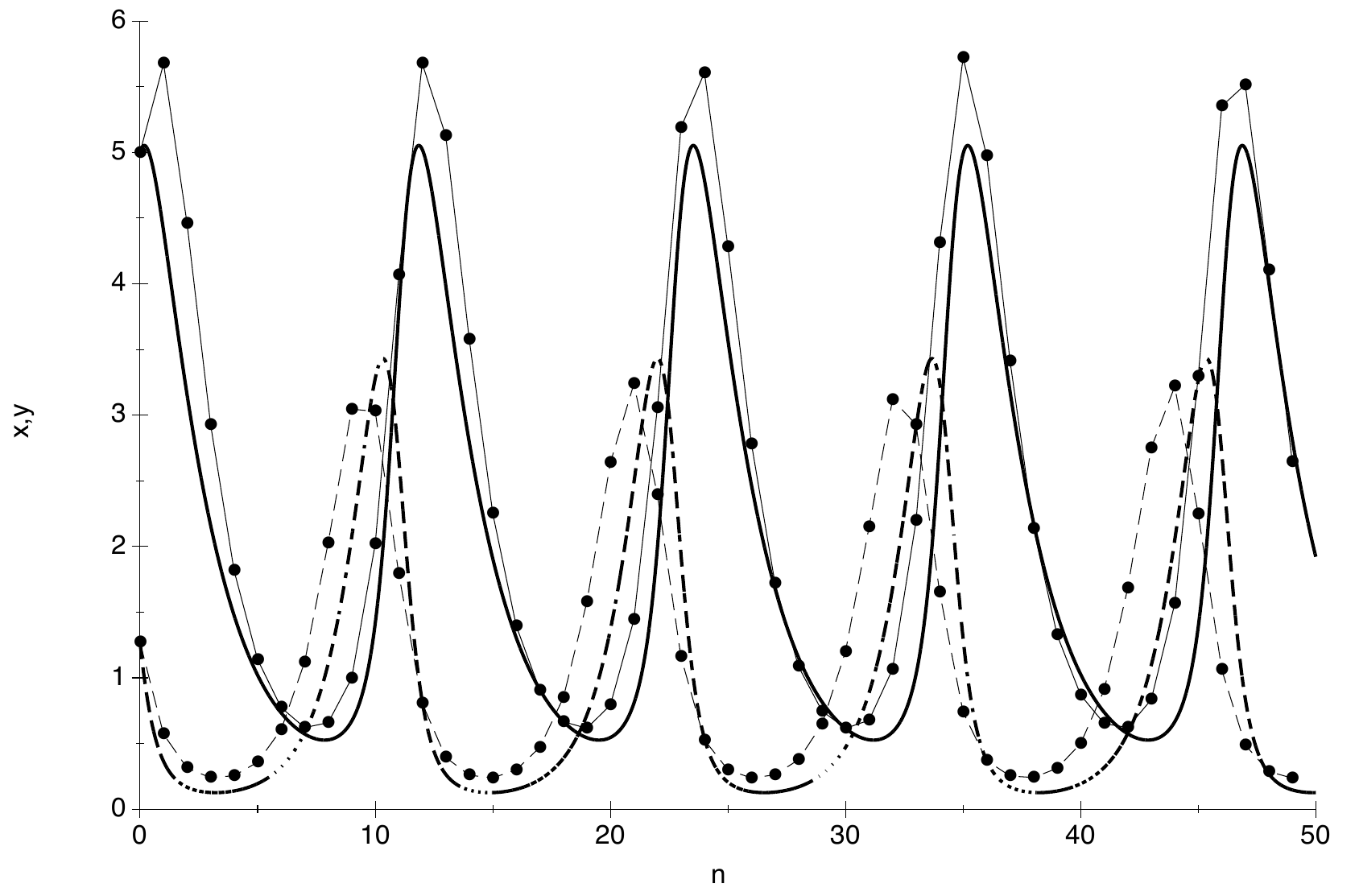}}
\smallskip{\bf Figure \six}. {\sl Numerical solution of the Lotka-Volterra equation, for $\delta=0.001$ (continuous and dashed lines for $x$ and $y$ respectively) and for $\delta=1$ (thin continuous and dashed lines, linking the points for $x$ and $y$ to enhance visibility)}.
\smallskip
It is remarkable that the evolution of the species still follows exactly the same pattern ---predator and prey are in quadrature, i.e. the maximum-minimum value of one coincides with the maximum absolute values of the derivative of the other--- even for values of the discretisation step that cover three orders of magnitude. In fact, the evolution keeps following the same pattern even when we increase the step for several more orders of magnitude. 

We know that in order to retrieve the initial differential system from its discrete analogue one must consider the limit where the discretisation step goes to zero. 
So one may wonder what the meaning could be of discretisation steps that attain values that have nothing to do with common discretisation practices. The answer to this question is to be found in a property of discrete systems for which positivity is guaranteed. Such systems possess another limit where the discretisation step can be taken to infinity, a limit known as the ultradiscrete limit [\refdef\ultra]. It allows, given a discrete system, to construct its generalised cellular automaton analogue. Let us show how this works in the case of the Lotka-Volterra system. Starting from equation (\vena) we introduce a small parameter $\epsilon$ through $\delta=e^{1/\epsilon}$. We rescale the variables $x$ and $y$ by $\delta$ and introduce new variables $X$ and $Y$ by setting $\delta x_n=e^{X_n/\epsilon}$ and $\delta y_n=e^{Y_n/\epsilon}$. Moreover we put $1+\delta\lambda=e^{L/\epsilon}$ and $1+\delta\mu=e^{M/\epsilon}$ where $L,M>0$. Taking the logarithm of (\vena) we find
$$X_{n+1}-X_n=L-\epsilon\log(1+e^{Y_n/\epsilon})\eqdaf\ttes$$
$$Y_{n+1}-Y_n=\epsilon\log(1+e^{X_{n+1}/\epsilon})-M.\eqno(\ttes{\rm b})$$
Next we take the limit $\epsilon\to+0$ and using the identity $\lim_{\epsilon\to+0}\epsilon\log(1+e^{Z/\epsilon})=\max(0,Z)$ we find
$$X_{n+1}=X_n+L-\max(0,Y_n)\eqdaf\tpen$$
$$Y_{n+1}=Y_n-M+\max(0,X_{n+1}).\eqno(\tpen{\rm b})$$
System (\tpen) defines a generalised cellular automaton in the sense that, if the initial conditions and the parameters are integers, the solution of the system remains integer at all time steps. This system is the ultradiscrete limit of the Lotka-Volterra system (\vena) and a simulation is given in Figure \figdef\seven.
\medskip
\centerline{\includegraphics[width=10 cm,keepaspectratio]{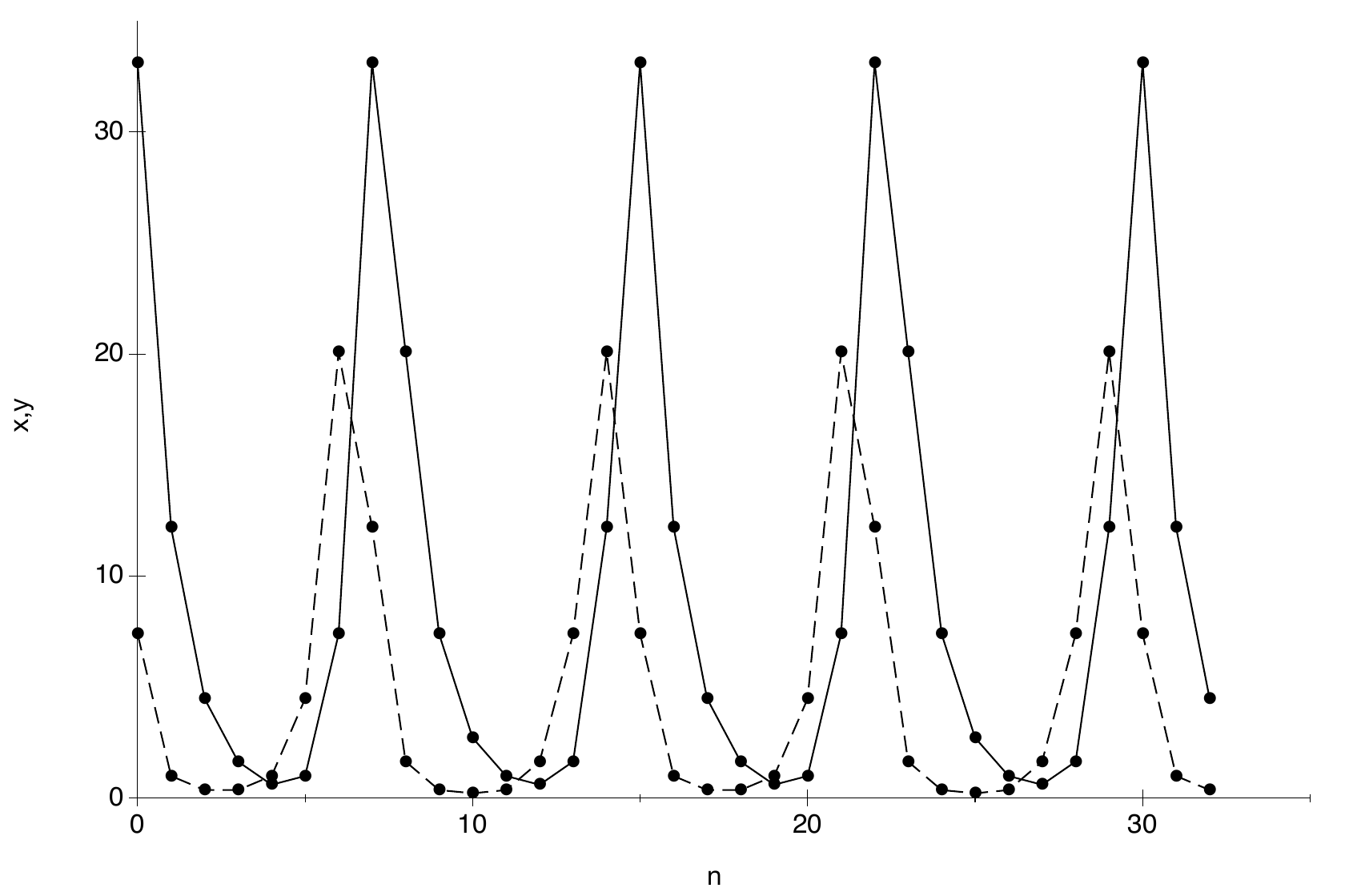}}
\smallskip{\bf Figure \seven}. {\sl Numerical solution of the ultradiscrete Lotka-Volterra equation}.
\smallskip
Note that in order to be able to visually compare the evolution obtained from (\tpen) to that given by (\vena), as shown in Figure 6, instead of $X$ and $Y$ we have plotted the quantities $e^{\sigma X}$ and $e^{\sigma Y}$, for some $\sigma>0$ chosen such that the magnitudes of the solutions in the plots are not too dissimilar.
It is remarkable that the evolution of the ultradiscrete system follows closely that of the discrete and the continuous one. This is not a unique occurence. As shown repeatedly [\refdef\epide, \refdef\oscil, \refdef\zultra] the dynamics of a cellular automaton constructed from a discrete system through the ultradiscretisation procedure follows closely that of its parent system, with perhaps, depending on the case, a few simplifications, but always keeping the essential features. 
\section{Integrable discretisations}
As we have seen in Section 3, one of Mickens' principles is that the discretisation must preserve the properties of the continuous system. This is of the utmost importance in the case of integrable systems. Already in Section 4, when we dealt with the quartic oscillator, we pointed out that the two discretisations we proposed did indeed possess a conserved quantity in the form of a discretised Hamiltonian. 

In this section we are going to concentrate on the discretisations of Painlev\'e equations, focusing in particular on the first two: P$_{\rm I}$ and P$_{\rm II}$. It is important to point out that in the discretisation of the Painlev\'e equations, in general, we do not aim from the outset at obtaining a correct {\sl explicit} dependence of the discrete system on the independent variable $n$, although the Painlev\'e equations themselves are of course non-autonomous equations. As we have explained in previous works of ours [\refdef\desot], this dependence is easily taken care of by the so-called deautonomisation process which, starting from an integrable autonomous discrete system,  consists in determining the possible non-autonomous forms of its coefficients through the application of some integrability criterion. In order to construct an autonomous discretisation we opt for a generalisation of Hirota's approach, allowing all forms of the proper homogeneity to exist, including local ones. Thus, for example, a term $\alpha x^2$ in a second-order equation can be discretised as $\alpha_1x_{n+1}x_{n-1}+\alpha_2x_n(x_{n+1}+x_{n-1})+\alpha_3x_n^2$. Moreover, in our approach, the parameter $\delta$ associated with the continuum limit at $\delta\to0$, need not necessarily be small.  As a consequence, the derivative term is replaced by a difference one with a coefficient which is on an equal footing with the remaining coefficients of the discrete equation. 
Since we will be focusing on Painlev\'e equations, the choice of autonomicity means that we deal, initially, with equations that have solutions given in terms of elliptic functions. In a purely discrete setting the QRT family of mappings [\qrt], the solution of which is given in terms of elliptic functions,  is a perfect starting point. The fact that we have a complete classification [\class] of canonical forms for the QRT mappings is of course a precious help. 

We start with the discretisation of the Painlev\'e I equation,
$$x''=x^2+ \lambda x + \mu(t),\eqdef\thex$$
where $\mu(t)$ is a linear function in the independent variable and the constant $\lambda$ should be set to zero in the canonical form of the Painlev\'e I equation. As explained above, since we will be focusing at first on the purely autonomous form of this equation, $\mu$ will be taken constant for now. Following the approach set out above, the second derivative is discretised to $x_{n+1}+x_{n-1}-2x_n$, while the quadratic term, $x^2$, is rendered as $ax_{n+1}x_{n-1}+bx_n(x_{n+1}+x_{n-1})+cx_n^2$. (Obviously no terms of the form $x_{n+1}^2, x_{n-1}^2$ can be present since we are looking for a well-defined reversible mapping).
Similarly, for the terms linear in $x$ we have $f(x_{n+1}+x_{n-1})+gx_n$. Finally our ansatz for a discrete form becomes
$$ax_{n+1}x_{n-1}+bx_n(x_{n+1}+x_{n-1})+cx_n^2+f(x_{n+1}+x_{n-1})+gx_n+h=0,\eqdef\thep$$
depending on 6 parameters.
Several constraints exist on these parameters in order to ensure that (\thep) does not become trivially linear. (They can be found in [\refdef\hirom]). 

We first examine the choice $a=0$ (putting $b=1$ by scaling). We have a first discrete system ($c=f=0$)
$$x_n(x_{n+1}+x_{n-1})+gx_n+h=0,\eqdef\toct$$
which has a conserved quantity
$$K=x_n^2x_{n+1}^2+gx_nx_{n+1}(x_n+x_{n+1})+(g^2+h)x_nx_{n+1}+gh(x_n+x_{n+1}).\eqdef\tenn$$
Upon deautonomisation we find that $h$ must be a linear function of $n$, and, of course, no invariant exists in this case.
Its canonical form is 
$$x_{n+1}+x_{n-1}=-g+{z_n\over x_n},\eqdef\tdek$$
where $z_n=\alpha n+\beta$, which is a well-known [\refdef\capel] discrete analogue of P$_{\rm I}$.
A second possibility exists, still for $a=0$, provided that $c=b$ (and $f=0$). The equation
$$x_n(x_{n+1}+x_n+x_{n-1})+gx_n+h=0,\eqdef\qena$$
has the conserved quantity
$$K=x_nx_{n+1}(x_n+x_{n+1})+gx_nx_{n+1}+h(x_n+x_{n+1}).\eqdef\qdyo$$
When $h$ is a linear function of $n$ the equation is again a discrete form of Painlev\'e I,
 $$x_{n+1}+x_n+x_{n-1}=-g+{z_n\over x_n}.\eqdef\qtri$$

Integrable cases exist also for $a\ne0$. 
The mapping
$$x_{n+1}x_{n-1}+gx_n+h=0,\eqdef\qtes$$
has a conserved quantity
$$K={x_nx_{n+1}(x_n+x_{n+1})-g(x_n^2+x_{n+1}^2)+(g^2-h)(x_n+x_{n+1})+gh\over x_nx_{n+1}}.\eqdef\qpen$$
The deautonomisation of (\qtes) leads to $h=\eta q^n$, which makes the equation a discrete, $q$-analogue, of P$_{\rm I}$.
The last case we are going to examine is
$$x_{n+1}x_{n-1}+bx_n(x_{n+1}+x_{n-1})+gx_n+h=0.\eqdef\qhex$$
We find that, when $b=1$, (\qhex) has the conserved quantity
$$K={x_n^2x_{n+1}^2+gx_nx_{n+1}(x_n+x_{n+1})+g^2x_nx_{n+1}-h(x_n^2+x_{n+1}^2)+h^2\over x_n+x_{n+1}}.\eqdef\qhep$$
The canonical form of (\qhex) after deautonomisation is
$$(x_{n+1}+x_n)(x_n+x_{n-1})=(x_n-z_n)^2-c^2\eqdef\qoct$$
where $z_n$ is again linear in $n$.

Several more integrable discretisation of the Painlev\'e I equation do exist, corresponding to more complicated canonical forms of the QRT mapping. Moreover, in our quest for discrete integable mappings we can find more integrable cases which are not analogues of a  Painlev\'e equation but mappings that can be linearised through some birational transformations. In this article we opted not to present these cases (although they are integrable, albeit definitely simpler than the discrete Painlev\'e).

Next we turn to the discretisation of Painlev\'e II
$$x''=x^3+\lambda(t) x+\mu\eqdef\qenn$$
where $\lambda(t)$ is now a linear function of the independent variable and $\mu$ a constant. Since we are in the presence of cubic nonlinearities we consider the discrete form
$$ax_{n+1}x_nx_{n-1}+bx_n^2(x_{n+1}+x_{n-1})+cx_n^3+f(x_{n+1}+x_{n-1})+gx_n+h=0,\eqdef\qdek$$
where again some constraints on the parameters are necessary for the mapping not to be trivial.

We start again with $a=0$ and setting $b=1$ find the equation
$$x_n^2(x_{n+1}+x_{n-1})+f(x_{n+1}+x_{n-1})+gx_n+h=0,\eqdef\pena$$
with conserved quantity
$$K=x_n^2x_{n+1}^2+gx_nx_{n+1}+f(x_n^2+x_{n+1}^2)+h(x_n+x_{n+1}).\eqdef\qhep$$
The deautonomisation of (\pena) shows that $g$ must be linear in $n$ with $f$ and $h$ constant. The canonical form of the discrete Painlev\'e equation obtained is
$$x_{n+1}+x_{n-1}={z_nx_n+h\over x_n^2-1}.\eqdef\pdyo$$
One utility of a discrete analogue of a Painlev\'e equation  that is unfortunately not stressed enough, is that it provides an excellent integrator of its continuous brethren, suitable for numerical calculations [\refdef\doriz]. In Figure \figdef\eight, we present a simulation of a solution of P$_{\rm II}$ that traverses multiple poles.
\medskip
\centerline{\includegraphics[width=10 cm,keepaspectratio]{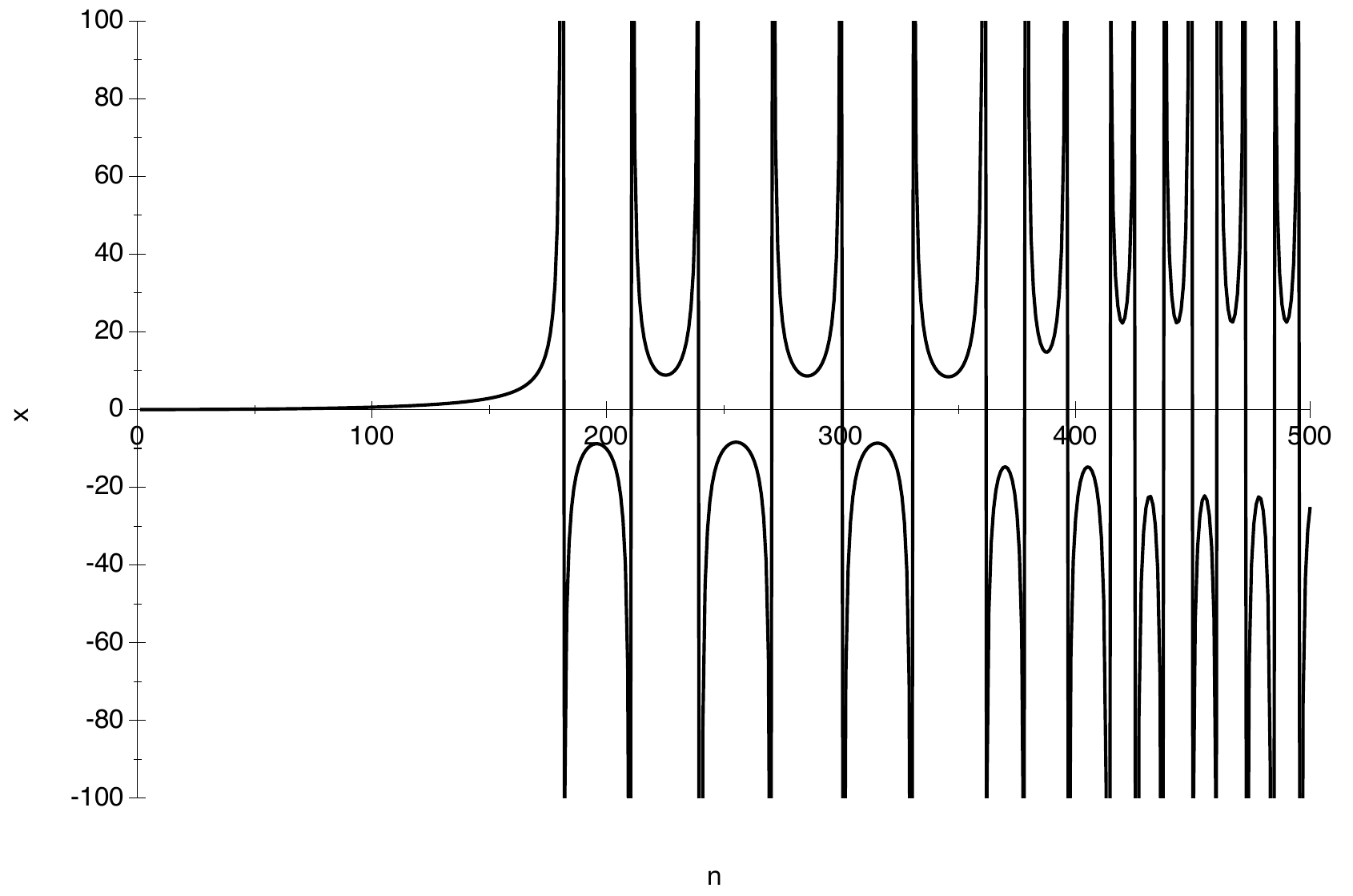}}
\smallskip{\bf Figure \eight}. {\sl A numerical simulation of the solution of Painlev\'e II displaying a succession of (simple) poles on the real axis}.
\smallskip

At this point we cannot resist the temptation to refer again to Kahan. In his classical paper [\kahan] he presents a solution of a Riccati equation displaying a series of poles and comments: {\sl ``What these results do not show is any sign that the computed trajectory had to pass through poles of the dependent  function"}. And he adds: {\sl ``Can any conventional method cope so simply with poles? None I know of''}. We feel that 
Kahan's remarks can be applied to the case of Painlev\'e II without changing a iota.

When $a\ne0$ we find, for example, the integrable mapping
$$x_{n+1}x_nx_{n-1}+gx_n+h=0,\eqdef\ptri$$
with conserved quantity 
$$K={x_nx_{n+1}(x_n+x_{n+1})-g(x_n+x_{n+1})-h\over x_nx_{n+1}}.\eqdef\ptes$$
With $g=\zeta q^n$ the equation in question becomes a $q$-discrete Painlev\'e equation, but it turns out that the continuum limit leads to P$_{\rm I}$.
A second integrable mapping can be obtained from 
$$ax_{n+1}x_nx_{n-1}+f(x_{n+1}+x_{n-1})+gx_n+h=0,\eqdef\ppen$$
provided $f=-a$. In this case (\ppen) has a conserved quantity of the form
$$K={x_n^2x_{n+1}^2-hx_nx_{n+1}(x_n+x_{n+1})+(g^2+h^2)x_nx_{n+1}-g(x_n^2+x_{n+1}^2)+gh(x_n+x_{n+1})\over x_nx_{n+1}-1}.\eqdef\phex$$
After deautonomisation, which gives $h_n=\eta q^n$, $g_n=\zeta q^{2n}$, the equation, written canonically as 
$$(x_{n+1}x_n-1)(x_nx_{n-1}-1)=-g_nx_n^2-h_nx_n+1\eqdef\phep$$
is indeed a $q$-discrete analogue of Painlev\'e II.

More examples of discrete analogues of the first two Painlev\'e equations can be given in the same vein. And, of course a similar approach can be implemented for the remaining equations of Painlev\'e. We shall not spend more time on this point in this paper, since our aim was merely to introduce discretisation methods that aim at preserving integrability. An exhaustive study of discrete Painlev\'e equations can be found in related publications of the authors. (See, for instance, [\refdef\cimpa] or [\refdef\physc]).
\section{Discussion}
In this paper we addressed the question of finding a suitable discretisation for a given differential system. We presented our work under the label ``unconventional'' discretisations, a term we borrowed from Kahan. What is meant by this is that a discretisation must be tailored to the system it is meant to discretise, precluding the use of general-purpose algorithms. It is our belief that a proper discretisation should give sensible results not only for a (very) small discretisation step, but must also lead to a physical behaviour of the solution even for a step close to 1.

One argument we insist upon in our approach is that of positivity. It is quite frequent in physics that the quantity one is dealing with, like the amplitude of a wave or the intensity of a field,  is positive. Examples of positive quantities in other domains like chemistry, biology or economics abound. When this is the case, it is imperative that the discretisation preserve the positivity of the solution. A practical rule as to how to achieve this was presented in the form of the ``aphorism'' in Section 5. The advantage of guaranteed positivity is that it makes it possible to take the ultradiscrete limit of a given discrete system (and by extension of the initial continuous one) transforming it into a generalised cellular automaton. 

One domain where ``unconventional'' discretisations play a dominant role is that of integrability. Following Mickens' principles, the discretisation of an integrable system must preserve the integrable character. This calls for special techniques and their development is an ongoing research. The present paper gave just a glimpse of these approaches, which go hand-in-hand with the elaboration of efficient discrete integrability detectors [\refdef\euler].

As discussed in the introduction, the mathematical description of physical phenomena proceeds through differential equations, based on the tacit assumption that space-time is continuous. This assumption was challenged already in ancient times. Zeno presented his dichotomy paradox claiming that in locomotion one must arrive at the half-way stage before he arrives at the goal, and so on, repeating ad infinitum the division by two. He concludied that, because an infinite number of steps are needed, the goal will never be reached. It is believed that Democritus countered the paradox by claiming the impossibility of infinite division of space and time (and, in fact of matter, as well). The discrete nature of space-time resolves the {\sl reductio ad absurdum} of Zeno's paradox.

If space-time were discrete this would entail that the equations describing the physical world be discrete.  So it is worth casting a look at what modern field theory could say on the matter. It is the common lore among physicists that the Planck length is the minimum physically meaningful distance. But this should {\sl not} be interpreted as meaning that the Planck length is like a pixel size of the physical world. In fact the Planck length is the length scale at which the Compton and the Schwarzschild lengths  of a hypothetical particle with a mass of the order of the Planck mass $\sqrt{\hbar c\over G}$ would coincide:
$$l_P=\sqrt{\hbar G\over c^3}\eqdef\poct$$ 
 where $\hbar$ is the Planck's constant divided by $2\pi$, $G$ the gravitational constant and $c$ the speed of light. It is of the order of $10^{-35}$ m and is thought to be the length scale at which quantum gravity becomes relevant. Lacking a consistent quantised theory for gravity we do not have any idea on what is happening for lengths smaller than the Planck length. Geometry may be affected by quantum uncertainties, it may even be non-commutative and space time may well be discrete. If this were true, tools like those presented in this article will provide a precious help in the description of the physical world.
 \paragraph{Acknowledgments.}
 RW  would like to thank the Japan Society for the Promotion of Science (JSPS) for financial support through the KAKENHI grant 23K22401. 
%
%

\begin{description}
 \item{[\verhu]} P.F. Verhulst, Notice sur la loi que la population suit dans son accroissement, Correspondance Math\'ematique et Physique (Ghent) 10 (1838) 113-121.
\item{[\logis]} E.N. Lorenz, The problem of deducing the climate from the governing equations, Tellus 16 (1964) 1-11.
\item{[\may]} R.M. May, Simple mathematical models with very complicated dynamics, Nature 261 (1976) 459-467.
\item{[\abram]} M. Abramowitz and I. Stegun, Handbook of Mathematical Functions, NBS, (1970), p. 896.
\item{[\skel]} J.G. Skellam, Random dispersal in theoretical populations, Biometrika 38 (1951) 196-218. 
\item{[\moris]} M. Morishita, The fitting of the logistic equation to the rate of increase of population density, Researches on Population Ecology 7 (1965) 52-55. 
\item{[\ronal]} R. E. Mickens, Exact solutions to a finite-difference model of a nonlinear reaction-advection equation: Implications for numerical analysis, Numer. Methods Partial Differ. Equ. 5, (1989) 313-325.
\item{[\micke]} R. E. Mickens, Advances in the Applications of Nonstandard Finite Difference Schemes, World Scientific, Singapore, (2005), p. 1.
\item{[\hirot]} R. Hirota, Nonlinear Partial Difference Equations: I. A Difference Analogue of the Korteweg-de Vries Equation, J. Phys. Soc. Japan 43, (1977) 1424-1433.
\item{[\kahan]} W. Kahan, {\sl Unconventional numerical methods for trajectory calculations},  unpublished lecture notes (1993).
\item{[\murat]} M. Murata, J. Satsuma, A. Ramani and B. Grammaticos, How to discretise differential systems in a systematic way, J. Phys. A 43 (2010) 315203, 15 pp.  
\item{[\potts]} R. Potts, Best difference approximation to Duffing's equation. J. Austral. Math. Soc., 23 (1982) 349-356.
\item{[\takah]} R. Hirota and D. Takahashi, {\sl Sabun to chorisan} (Kyoritsu Shuppan, Tokyo, 2003) [{\sl Discrete and ultradiscrete systems}, in Japanese].
\item{[\suris]} Yu. Suris, Integrable mappings of the standard type, Funct. Anal. Appl., 23 (1989) 74-76.
\item{[\class]} A. Ramani, S. Carstea, B. Grammaticos and Y. Ohta, On the autonomous limit of discrete Painlev\'e equations, Physica A  305 (2002) 437-444.
\item{[\qrt]} G.R.W. Quispel, J.A.G. Roberts and C.J. Thompson, Integrable mappings and soliton equations, Physica D34 (1989) 183-192.
\item{[\polar]} E. Celledoni, R. McLachlan, D. McLaren, B. Owren and G.R.W. Quispel, Discretizationof polynomial vector ﬁelds by polarization, Proc. Royal Soc. A. 471 (2015), no.2184, 20150390, 10 pp.
\item{[\yuris]} Yu. Suris, A new approach to integrals of discretizations by polarization, Op. Comm. Nonlin. Math. Phys. Special Issue 1, (2024) 1-8.
\item{[\handy]} B. Grammaticos, R. Willox and J. Satsuma, Revisiting the Human and Nature Dynamics model, Reg. Chao. Dyn. 25 (2020) 178-198. 
\item{[\delay]} A. Ramani, B. Grammaticos, J. Satsuma and R. Willox, Discretisation induced delays and their role in the dynamics, J. Phys. A 41 (2008) 205204, 11 pp.
\item{[\ultra]} T. Tokihiro, D. Takahashi, J. Matsukidaira and J. Satsuma, From Soliton Equations to Integrable Cellular Automata through a Limiting Procedure, Phys. Rev. Lett. 76 (1996) 3247-3250.
\item{[\epide]} R. Willox, B. Grammaticos, A.S. Carstea and A. Ramani, Epidemic dynamics: discrete-time and cellular automaton models, Physica A 328 (2003) 13-22.
\item{[\oscil]} A. Ramani, A.S. Carstea, R. Willox and B. Grammaticos, Oscillating epidemics: a discrete-time model, Physica A 333 (2004) 278-292.
\item{[\zultra]} R. Willox, Y. Nakata, J. Satsuma, A. Ramani and B. Grammaticos, Solving the ultradiscrete KdV equation, J. Phys. A 43 (2010) 482003, 7 pp.
\item{[\desot]} T. Mase, R. Willox, B. Grammaticos and A. Ramani, Deautonomisation by singularity confinement: an algebro-geometric justification, Proc. Roy. Soc. A 471 (2015) 20140956, 20 pp.
\item{[\hirom]} B. Grammaticos, A. Ramani, J. Satsuma and R. Willox, Discretising the Painlev\'e equations \`a la Hirota-Mickens, J. Math. Phys. 53 (2012) 023506, 24 pp.
\item{[\capel]} A. Ramani and B. Grammaticos, Discrete Painlev\'e equations: coalescences, limits and degeneracies, Physica A 228 (1996) 160-171.
\item{[\doriz]} B. Grammaticos and B. Dorizzi, Integrable discrete systems and numerical integrators, J. Math. Comp. in Sim. 37 (1994) 341-352.
\item{[\cimpa]} B. Grammaticos and A. Ramani, Discrete Painlev\'e equations: a review, Springer LNP 644 (2004) 245-321.
\item{[\physc]} B. Grammaticos and A. Ramani, Discrete Painlev\'e equations: an integrability paradigm, Phys. Scr. 89 (2014) 038002, 13 pp.
\item{[\euler]} B. Grammaticos, A. Ramani, R. Willox and  T. Mase, in {\sl Nonlinear Systems and Their Remarkable Mathematical Structures,} N. Euler (ed), CRC Press, Boca Raton FL, (2018), p. 44-74.
\end{description}

\end{document}